\documentclass[12pt]{amsart}

\usepackage[a4paper,left=20mm,right=20mm,top=25mm,headheight=5mm,headsep=5mm,bottom=25mm,footskip=10mm]{geometry}
\usepackage[T1]{fontenc} 
\usepackage[active]{srcltx}
\usepackage{amsthm,amsmath,amssymb,amsfonts,amscd}
 \usepackage{bbm}
\usepackage{hyperref}
\hypersetup{colorlinks=true,allcolors=blue}
\usepackage{graphicx}
\usepackage{enumerate}
\usepackage{enumitem}
\usepackage[dvipsnames]{xcolor}
\usepackage{physics}
\usepackage{pmboxdraw}
\usepackage{listings}
\usepackage{tikz}
\usepackage{mathtools}
\usepackage{float}
\usepackage{xifthen}
\usepackage{bm}
\usepackage{subfig}
\usepackage{inconsolata}
\usepackage{centernot}
\usepackage[framemethod=tikz]{mdframed}

\usepackage{subfig}
\usepackage{comment}

\newtheorem{theorem}{Theorem}[section]

\numberwithin{equation}{subsection}
\numberwithin{figure}{subsection}
\numberwithin{theorem}{subsection}

\newtheorem{lemma}[theorem]{Lemma}
\newtheorem{corollary}[theorem]{Corollary}

\theoremstyle{definition}

\newtheorem{example}[theorem]{Example}

\theoremstyle{remark}

\newtheorem{remark}[theorem]{Remark}

\theoremstyle{plain}
\newmdtheoremenv[%
   backgroundcolor=cyan!10,%
   outerlinecolor=black,%
   innerbottommargin = \topskip,%
   roundcorner=4]
{QUESTION}{Question}

\newcommand{\N}{\mathbb{N}}
\newcommand{\Z}{\mathbb{Z}}

\newcommand{\C}{\mathbb{C}}

\DeclareMathOperator{\Gaps}{Gaps}
\DeclareMathOperator{\ord}{ord}

\definecolor{blue-violet}{rgb}{0.54, 0.17, 0.89}
\definecolor{Blue}{rgb}{0.01, 0.28, 1.0}
\definecolor{gGreen}{rgb}{0.2, 0.8, 0.2}
\definecolor{Green}{rgb}{0.04, 0.85, 0.32}

\begin{document}

\title[Analytic information encoded in the Bernstein--Sato polynomial]{Analytic information encoded in the Bernstein--Sato polynomial of genus-one plane branches}



\author[M. Alberich-Carrami\~nana]{Maria Alberich-Carrami\~nana}
\address{
    Departament de Matem\`atiques  and  Institut de Matem\`atiques de la UPC-BarcelonaTech (IMTech)\\  Universitat Polit\`ecnica de Catalunya \\ Av.~Diagonal 647, Barcelona 08028; and Institut de Rob\`otica i Inform\`atica Industrial\\ CSIC-UPC \\ Llorens i Artigues 4-6, Barcelona 08028, Spain.
    } 
\email{Maria.Alberich@upc.edu}

\author[M. Morella]{Maria Morella Gim\'enez}
\address{
    Facultat de Matem\`atiques  i Estadística\\  Universitat Polit\`ecnica de Catalunya \\ C.~Pau Gargallo 14, Barcelona 08028.
    } 
\email{mamogi2001@gmail.com}

\thanks{
MAC is partially supported by the Spanish State Agency project PID2023-146936NB-I00 funded by MICIU/AEI/10.13039/501100011033 and, by ERDF/EU, and by the GEOMVAP 2021-SGR-00603 AGAUR project.}


\begin{abstract}
We give an explicit formula recovering the Zariski invariant of a plane branch with one Puiseux pair from its Bernstein–Sato polynomial. Thus the Bernstein–Sato stratification refines the Zariski-invariant stratification. We determine the generic roots in each stratum from the value semigroup and the Zariski invariant, refining Yano’s description for a fixed topological class. We also compute the maximal spectral–exponent discrepancy and relate it to the dimension of the generic Bernstein–Sato locus and the minimum Tjurina number. 
%
\end{abstract}

\maketitle

\section{Introduction}

One of the central challenges in singularity theory is determining what analytic information is encoded in a given analytic invariant, especially which information can be retrieved.
Among discrete analytic invariants, the Bernstein--Sato polynomial of a germ of a holomorphic function $f \in \mathcal{O}_n = \mathbb{C}\{x_1, \dots , x_n \}$ is distinguished for its connections to D-module and Hodge theory, though the analytic information it encodes is still only partially understood. 
Recall that, the 
\textit{Bernstein--Sato polynomial} $b_{f}$ of $f$ is the monic polynomial of the smallest degree in $\mathbb{C}[s]$ 
fulfilling the functional equation
\[
P(s) \cdot f^{s+1} = b_{f}(s) f^s,
\]
where $P(s)$ is a differential operator with holomorphic coefficients in $x_1, \dots, x_n$ and polynomial dependence in $s$ (see \cite{bernstein1972analytic,bjork1974dimensions,sato1970theory}).
%
Since \(-1\) is always a root of \(b_f(s)\), we denote by $\mathcal R_f$
the multiset  of roots (counted with their algebraic multiplicities) of the reduced Bernstein--Sato polynomial \(\widetilde b_f(s)=b_f(s)/(s+1)\). This set 
satisfies $\mathcal{R}_f \subset \mathbb{Q}_{<0}$ 
and its cardinal is at most the Milnor number $ \mu _f$ of $f$ \cite{kashiwara-b-functions, malgrange1975polynome}. 


Malgrange \cite{malgrange1975polynome} characterized \(\widetilde b_f(s)\)
as the minimal polynomial of an endomorphism induced by the Gauss--Manin
connection on a \(\mu_f\)-dimensional vector space related to the Brieskorn
lattice of \(f\). Let \(\mathcal B_f\) be the multiset of roots of the
characteristic polynomial of this endomorphism, called the
\emph{\(b\)-exponents}. In particular, \(\#\mathcal B_f=\mu_f\).
If the eigenvalues of the monodromy of \(f\) are pairwise distinct,
as will be the case throughout this paper, the characteristic polynomial
equals \(\widetilde b_f(s)\). Hence \(\mathcal B_f=\mathcal R_f\).


The Bernstein--Sato polynomial is not, in general, a topological
invariant. For plane \textit{branches}, i.e. germs of curves defined by an
irreducible \(f\in\mathbb{C}\{x,y\}\), some of its roots remain constant
within a fixed topological class (called \emph{topological roots}), whereas others may change under \(\mu_f\)-constant deformations, changing in a jumping way by a shift of $\pm 1$. 
Understanding the analytic 
nature of these variations, or the lack thereof, has been the focus of significant recent research. 
%

Blanco identified topological roots of geometric origin for arbitrary plane branches \cite{blanco2024topological}, and Yano's conjecture, proved by Blanco in full generality and by Cassou-Nogu\`es in the \textit{genus-one} (one Puiseux-pair) case, predicts the \(b\)-exponents occurring generically in a fixed topological class \cite{yano1982exponents,blanco2021yano,cassou1987etude}.
Recent works by Senovilla \cite{senovilla2024differential} and Guimarães–Hernandes \cite{guimaraes2025bernstein} characterize, in the genus-one case, some Bernstein--Sato roots from the semimodule of values of K\"ahler differentials, which is an analytical invariant of the branch. 


In this work we provide a new insight in the more interesting and less explored direction: what the Bernstein-Sato roots reveals about the analytical class of the branch, namely about the semimodule of values of K\"ahler differentials (which equals the semimodule of differential values). Let \(\gamma\) be a plane branch. Its semigroup of values \(\Gamma_\gamma
\) is a complete invariant of the topological class (meaning that \(\Gamma_\gamma
\) characterizes the topological class), while the semimodule of differential values \(\Lambda_\gamma
\) is the basic discrete invariant used in the analytic classification of branches; see, for instance, \cite{hefez2011analytic}. 

The non-topological part \(\Lambda_\gamma\setminus \Gamma_\gamma\)
measures analytic information not visible in the embedded topology of $\gamma$. For a branch $\gamma$ with one Puiseux-pair $\frac{m}{n}$, the first such value is
\[
   \lambda_1:=\min(\Lambda_\gamma\setminus \Gamma_\gamma),
\]
with the convention \(\lambda_1=\infty\) if \(\Lambda_\gamma\setminus\Gamma_\gamma=\emptyset\). It is related to the Zariski invariant  \(j_1\) \cite{zariski1966maximumtorsion} by \cite{hefez2011analytic,casas2023polar}
\[
   \lambda_1=j_1+n+m.
\]

The first main result of this paper is that the Bernstein--Sato polynomial determines \(\lambda_1\), and hence \(j_1\), for all genus-one plane branches (Theorem~\ref{MainThm}). This is a recovery result in the opposite direction to the recent theorems of Senovilla and Guimar\~aes--Hernandes, which show that  differential values and spectral numbers determine some Bernstein--Sato roots \cite{senovilla2024differential,guimaraes2025bernstein}. The recovery is optimal: the Bernstein--Sato polynomial does not determine further non-topological differential values in general.

The recovery result in Theorem~\ref{MainThm} is particularly noteworthy because no other discrete analytic invariant, aside from the differential values, has been observed to exhibit such a property.
For instance, the Tjurina number $\tau_\gamma$, which is defined in terms of the number of gaps of the semimodule of differentials, namely $\tau_\gamma = \mu_\gamma - \# \Lambda_\gamma \setminus \Gamma_\gamma$, 
does not provide sufficient analytic information to determine the Zariski invariant, as noted by Peraire in~\cite{peraire1998moduli}. Similarly, the topological class of the generic polars of a genus-one plane branch does not determine the Zariski invariant, as demonstrated by Casas-Alvero in~\cite{casas2023polar}.

We use Theorem~\ref{MainThm} to establish the main contribution of this paper in Theorem~\ref{thm:generic-roots}: a description of the generic Bernstein--Sato roots in each Zariski-invariant stratum of a fixed topological class, that is, within the moduli space \(\mathcal M_\Gamma\) of analytic classes of branches with fixed semigroup \(\Gamma=\langle n,m\rangle\). Since the Bernstein--Sato polynomial determines the Zariski invariant \(j_1\), the stratification by Bernstein--Sato roots refines the stratification by \(j_1\) studied by Peraire \cite{peraire1998moduli}. We prove that, in each stratum with fixed \(j_1\), the set of Bernstein--Sato roots attained by generic branches can be predicted from \(\Gamma\) and \(j_1\) refining Yano’s description for a fixed
topological class. 
As consequences, we obtain a closed formula for the maximal discrepancy between spectral numbers and exponents (associated with the reduced Bernstein--Sato roots) in that stratum and 
we identify the dimension of the generic Bernstein--Sato substratum with the difference between this maximal discrepancy and \(\mu-\tau_{\min}(\Gamma,j_1)\), where \(\tau_{\min}(\Gamma,j_1)\) is the minimal Tjurina number among branches in that stratum, attained by the generic ones.

\section{Preliminaries for genus-one branches}

Let \(\gamma:f=0\) be a genus-one plane branch, i.e. with one Puiseux-pair $\frac{m}{n}$. We write
\[
   \Gamma=\Gamma_\gamma =\Gamma_f=\langle n,m\rangle,
   \qquad 1<n<m,
   \qquad \gcd(n,m)=1.
\]
Thus the Milnor number of $\gamma$ is \(\mu=(n-1)(m-1)\), and the conductor of \(\Gamma\) is also \(\mu\). We use the notation
\[
   \Gaps(\Gamma):=\N\setminus \Gamma.
\]

Let \(\varphi(t)=(x(t),y(t))\) be a Puiseux parametrization of \(\gamma\). For \(g\in\C\{x,y\}\), set
\[
   \ord_\gamma(g):=\ord_t(g\circ\varphi).
\]
For a holomorphic one-form \(\omega=A(x,y)dx+B(x,y)dy\), its value along \(\gamma\) is
\[
   \ord_\gamma(\omega):=\ord_t(\varphi^*\omega)+1.
\]
The semimodule of differential values (equiv. the values of the Kähler differentials) is
\[
   \Lambda_f:=\{\ord_\gamma(\omega):\omega\in\Omega^1_{\C^2,0}\}.
\]
It is a \(\Gamma\)-semimodule, and \(\Gamma\setminus\{0\}\subseteq\Lambda_f\). The elements of \(\Lambda_f\setminus\Gamma\) are the \textit{non-topological differential values}.

Every branch in the topological class \(\Gamma=\langle n,m\rangle\) admits, after analytic changes of coordinates, a semi-quasihomogeneous equation of the form
\begin{equation}
   f=x^m+y^n+
     \sum_{j\in J} z_j x^{p_{1,j}}y^{p_{2,j}},
   \qquad z_j\in\C,
   \label{eq:normal-form}
\end{equation}
where the set of \textit{standing points} \cite[\S~2]{casas2023polar} is
\begin{equation}
   P=P_{n,m}:=\bigl\{(p_1,p_2)\in\Z_{\geq0}^2:
      0\leq p_1\leq m-2,
      0\leq p_2\leq n-2,
      np_1+mp_2>nm
   \bigr\},
   \label{eq:P}
\end{equation}
and the set indexing the \textit{standing coefficients} is
\begin{equation}
   J=J_{n,m}:=\bigl\{j=np_1+mp_2-nm:(p_1,p_2)\in P\bigr\}.
   \label{eq:J}
\end{equation}
For each \(j\in J\), the corresponding point of \(P\) is denoted by \((p_{1,j},p_{2,j})\). The Zariski invariant \cite{zariski1966maximumtorsion} is
\begin{equation}
   j_1=j_1(f):=\min\{j\in J:z_j\neq0\}.
   \label{eq:zariski-invariant}
\end{equation}
For the quasihomogeneous branch \(x^m+y^n=0\), all \(z_j\) vanish and we set \(j_1=\infty\). We also set
\begin{equation}
   \lambda(j):=j+n+m,
   \qquad s_j:=\frac{\lambda(j)}{nm}
   \qquad (j\in J),
   \label{eq:lambda-sj}
\end{equation}
and, when \(j_1<\infty\), by the descriptions of the Zariski invariant and differential values in \cite{hefez2011analytic,casas2023polar},
\begin{equation}
   \lambda_1:= \min(\Lambda_\gamma\setminus \Gamma)= j_1+n+m.
   \label{eq:lambda1}
\end{equation}
If \(j_1=\infty\), we set \(\lambda_1=\infty\).

The set used by Cassou-Nogu\`es in \cite{cassou1987etude} is the image of \(P\) under the central symmetry with respect to \(\bigl(\frac{m-1}{2},\frac{n-1}{2}\bigr)\):
\begin{equation}
   M=M_{n,m}:=\bigl\{(m-p_1-1,n-p_2-1):(p_1,p_2)\in P\bigr\}.
   \label{eq:M}
\end{equation}

We shall compare three root sets attached to the topological class determined by $\Gamma$. First, the topological roots of geometric origin in the genus-one case \cite{blanco2024topological}
 are
\begin{equation}
   \mathcal T_\Gamma:=
   \left\{
      -\frac{n+m+\nu}{nm}:
      \nu\in\Gamma,
      0\leq \nu<nm,
      n\frac{n+m+\nu}{nm}\notin\Z,
      m\frac{n+m+\nu}{nm}\notin\Z
   \right\}.
   \label{eq:T-roots}
\end{equation}
Second, the Yano set for the topological class \cite{blanco2021yano} is
\begin{equation}
   \mathcal Y_\Gamma:=
   \left\{
      -\frac{n+m+\nu}{nm}:
      0\leq \nu<nm,
      n\frac{n+m+\nu}{nm}\notin\Z,
      m\frac{n+m+\nu}{nm}\notin\Z
   \right\}.
   \label{eq:Y-roots}
\end{equation}
Third, the reduced roots of the quasihomogeneous branch \cite{cassou1987etude} are
\begin{equation}
   \mathcal Q_\Gamma:=
   \left\{
      -\frac{n(\alpha+1)+m(\beta+1)}{nm}:
      0\leq\alpha\leq m-2,
      0\leq\beta\leq n-2
   \right\}.
   \label{eq:Q-roots}
\end{equation}
The root \(-1\) of the full Bernstein--Sato polynomial is not included in these reduced root sets.

\begin{lemma}[Combinatorics of the three root sets]
\label{lem:combinatorics-root-sets}
With the notation above, the following equalities hold:
\begin{align}
   \mathcal Y_\Gamma\setminus\mathcal T_\Gamma
      &=\left\{-\frac{\lambda(j)}{nm}:j\in J\right\},
      \label{eq:Y-minus-T}\\
   \mathcal Q_\Gamma\setminus\mathcal T_\Gamma
      &=\left\{-\frac{\lambda(j)}{nm}-1:j\in J\right\}.
      \label{eq:Q-minus-T}
\end{align}
Moreover, the translation \(j\mapsto\lambda(j)=j+n+m\) induces a
bijection
\begin{equation}
   \lambda:J\xrightarrow{\sim}
   \{\ell\in\Gaps(\Gamma):\ell\geq n+m\}.
   \label{eq:J-large-gaps}
\end{equation}
The three sets
\[
   \mathcal T_\Gamma,\qquad
   \mathcal Y_\Gamma\setminus\mathcal T_\Gamma,\qquad
   \mathcal Q_\Gamma\setminus\mathcal T_\Gamma
\]
are pairwise disjoint. 
In particular,
\begin{equation}
   \#\mathcal T_\Gamma+\#J=(n-1)(m-1)=\mu.
   \label{eq:T-plus-J-mu}
\end{equation}
\end{lemma}

\begin{proof}
A positive integer \(\nu<nm\) is a gap of \(\Gamma=\langle n,m\rangle\) if and only if it can be written uniquely as
\[
   \nu=nm-an-bm,
   \qquad 1\leq a\leq m-1,
   \qquad 1\leq b\leq n-1.
\]
For such a gap, the condition that
\(-\frac{n+m+\nu}{nm}\) belongs to \(\mathcal Y_\Gamma\) is equivalent to asking that neither \(n+m+\nu\) be divisible by \(n\) nor by \(m\). Since \(\gcd(n,m)=1\), this is equivalent to \(a\geq2\) and \(b\geq2\). Setting \(p_1=m-a\) and \(p_2=n-b\) gives precisely a point of \(P\), and
\[
   \nu=np_1+mp_2-nm \in J.
\]
This proves \eqref{eq:Y-minus-T}. For \eqref{eq:Q-minus-T}, if \(j=np_{1,j}+mp_{2,j}-nm \in J\), then
\[
   -\frac{\lambda(j)}{nm}-1
   =-\frac{n(p_{1,j}+1)+m(p_{2,j}+1)}{nm},
\]
which is the quasihomogeneous root indexed by \((p_{1,j},p_{2,j})\).
Conversely, let a root of \(\mathcal Q_\Gamma\) be indexed by
\((p_1,p_2)\) in the rectangle from \eqref{eq:Q-roots}, and put
\(w=np_1+mp_2\). The equality \(w=nm\) is impossible: reduction modulo
\(m\) would give \(m\mid p_1\), hence \(p_1=0\), and then \(p_2=n\),
contrary to \(p_2\leq n-2\). If \(w<nm\), then
\[
   -\frac{n(p_1+1)+m(p_2+1)}{nm}
   =-\frac{n+m+w}{nm}
   \in\mathcal T_\Gamma.
\]
Indeed, \(w\in\Gamma\), and the two required nondivisibility conditions
follow from \(1\leq p_1+1\leq m-1\) and
\(1\leq p_2+1\leq n-1\). If \(w>nm\), then
\((p_1,p_2)\in P\), and \(j=w-nm\in J\) gives the displayed root
\(-\lambda(j)/(nm)-1\), which cannot belong to
\(\mathcal T_\Gamma\), since its parameter \(w\) is greater than \(nm\).
This proves \eqref{eq:Q-minus-T}.

We next prove the bijection \eqref{eq:J-large-gaps}. Given \(j\in J\),
write \(j=np_1+mp_2-nm>0\) with \((p_1,p_2)\in P\), and set
\[
   a=m-p_1-1,\qquad b=n-p_2-1.
\]
The bounds defining \(P\) give
\(1\leq a\leq m-1\) and \(1\leq b\leq n-1\), and
\[
   \lambda(j)=j+n+m=nm-na-mb>n+m.
\]
Thus \(\lambda(j)\) is a gap by the gap parametrization above.
Conversely, let \(\ell\in\Gaps(\Gamma)\) with \(\ell\geq n+m\).
Since \(n+m\in\Gamma\), in fact \(\ell>n+m\). Write
\[
   \ell=nm-na-mb,\qquad
   1\leq a\leq m-1,\quad 1\leq b\leq n-1,
\]
and put \(p_1=m-a-1\), \(p_2=n-b-1\). Then
\[
   0\leq p_1\leq m-2,\qquad 0\leq p_2\leq n-2,
   \qquad np_1+mp_2=nm+\ell-n-m>nm.
\]
Hence \((p_1,p_2)\in P\), and
\(j=np_1+mp_2-nm=\ell-n-m\) belongs to \(J\) and satisfies
\(\lambda(j)=\ell\). This proves surjectivity; injectivity follows
because \(\lambda\) is a translation.

It remains only to check the additional disjointness assertion. Since
\(\lambda(j)\) is a gap and the Frobenius number of \(\Gamma\) is
\(nm-n-m\), one has
\(
   0<\lambda(j)<nm.
\)
Therefore
\[
   \mathcal Y_\Gamma\setminus\mathcal T_\Gamma\subset(-1,0),
   \qquad
   \mathcal Q_\Gamma\setminus\mathcal T_\Gamma\subset(-2,-1),
\]
so these two sets are disjoint. Their intersections with
\(\mathcal T_\Gamma\) are empty by definition, proving that the three
sets are pairwise disjoint.
Finally, as \(\nu\) runs through \(0,\ldots,nm-1\), the integers
\(n+m+\nu\) run through a complete residue system modulo \(nm\).
Exactly \((n-1)(m-1)\) of these are divisible by neither \(n\) nor \(m\).
Hence \(\#\mathcal Y_\Gamma=(n-1)(m-1)=\mu\), and
\eqref{eq:Y-minus-T} yields \eqref{eq:T-plus-J-mu}.
\end{proof}

\subsection{Cassou-Nogu\`es' residue criterion}

For \(j\in J\) and \((a,b)\in M\), Cassou-Nogu\`es associates to the local zeta integral
\[
   I(f,(a,b))(s)=\iint f(x,y)^s x^{a-1}y^{b-1}\,dx\,dy
\]
a residue at \(s_j=\frac{j+n+m}{nm}\). In the normal form \eqref{eq:normal-form}, this residue is the following polynomial expression in the coefficients \(z_\ell\), whose weights involve the Euler gamma function $\Gamma (\cdot)$:
\begin{equation}
\begin{aligned}
   \operatorname{Res}_{f}(a,b)(s_j)
   &=\frac{\Gamma(s_j)^{-1}}{nm}
   \sum_{\substack{(\delta_\ell)_{\ell\in J}\in(\Z_{\geq0})^{\#J}\\
      \sum_{\ell\in J}\delta_\ell\ell=k_j(a,b)}}
   (-1)^{\sum_\ell\delta_\ell}
   \Gamma\!\left(\frac{\sum_\ell\delta_\ell p_{1,\ell}+a}{m}\right)
   \Gamma\!\left(\frac{\sum_\ell\delta_\ell p_{2,\ell}+b}{n}\right)
   \prod_{\ell\in J}\frac{z_\ell^{\delta_\ell}}{\delta_\ell!},
\end{aligned}
\label{eq:cassou-residue}
\end{equation}
where
\begin{equation}
   k_j(a,b):=s_jnm-na-mb=j+n(1-a)+m(1-b).
   \label{eq:k-j-ab}
\end{equation}
If no family \((\delta_\ell)\) satisfies the displayed condition, the sum is understood to be zero.


By Lemma \ref{lem:combinatorics-root-sets}, Cassou-Noguès’ residue criterion admits the following sharper formulation, in terms of a disjoint decomposition of sets of Bernstein–Sato roots.

\begin{theorem}[Cassou-Nogu\`es \cite{cassou1987etude}]
\label{thm:cassou-nogues}
Let \(f\) define a semi-quasihomogeneous branch in the topological class \(\Gamma=\langle n,m\rangle\), written as in \eqref{eq:normal-form}. Then the reduced Bernstein--Sato roots are described by the disjoint union
\begin{equation}
\begin{aligned}
   \mathcal R_f
   ={}&\mathcal T_\Gamma\\
   &\sqcup
   \left\{-s_j:j\in J,\ \operatorname{Res}_{f}(a,b)(s_j)\neq0
      \text{ for some }(a,b)\in M\right\}\\
   &\sqcup
   \left\{-s_j-1:j\in J,\ \operatorname{Res}_{f}(a,b)(s_j)=0
      \text{ for every }(a,b)\in M\right\}.
\end{aligned}
\label{eq:cassou-root-decomposition}
\end{equation}
Equivalently, for each \(j\in J\), the pair \(\{-s_j,-s_j-1\}\) contributes exactly one reduced Bernstein--Sato root: \(-s_j\) occurs if and only if at least one of the residues \(\operatorname{Res}_{f}(a,b)(s_j)\), \((a,b)\in M\), is nonzero; otherwise \(-s_j-1\) occurs.
\end{theorem}

The following result identifies the first non-topological differential value with the first nonzero residue of type \((1,1)\).

\begin{theorem}[Senovilla \cite{senovilla2024differential}]
\label{thm:senovilla}
Let \(f\) be as in \eqref{eq:normal-form}. For \(j\in J\), the following are equivalent:
\begin{enumerate}[label=\textup{(\roman*)}]
   \item \(\lambda(j)=j+n+m\) is the first non-topological differential value \(\lambda_1=\min(\Lambda_f\setminus\Gamma)\);
   \item \(j\) is the Zariski invariant \(j_1=\min\{i\in J:z_i\neq0\}\);
   \item \(\operatorname{Res}_{f}(1,1)(s_\ell)=0\) for all \(\ell<j\), \(\ell\in J\), and
   \[
      \operatorname{Res}_{f}(1,1)(s_j)\neq0.
   \]
\end{enumerate}
Moreover, if \(\lambda_1<\infty\), then for every \(\lambda=\lambda_1+na+mb\in(\lambda_1+\Gamma)\setminus\Gamma\) one has
\[
   \operatorname{Res}_{f}(a+1,b+1)\!\left(\frac{\lambda}{nm}\right)\neq0.
\]
\end{theorem}

Combining Theorems~\ref{thm:cassou-nogues} and~\ref{thm:senovilla}, Senovilla obtains that every 
element  of \((\lambda_1+\Gamma)\setminus\Gamma\) gives a Bernstein--Sato root.
%
{Guimar\~aes and Hernandes generalize this result to any non-topological differential value and prove the following two root inclusions by comparing filtrations on the Milnor algebra and the saturated Brieskorn lattice \cite[Theorem~4.1 and Proposition~4.4]{guimaraes2025bernstein}. The spectrum \(\mathcal S_f\), used in the second inclusion, is introduced in Section~\ref{sec:spectral-exponent-discrepancy}.}

\begin{theorem}[Guimar\~aes--Hernandes \cite{guimaraes2025bernstein}]
\label{thm:guimaraes-hernandes}
Let \(f\) define a semi-quasihomogeneous branch with semigroup
\(\Gamma=\langle n,m\rangle\), and let \(\mathcal S_f\) be its spectrum
in the \((-1,1)\) normalization. Then:
\begin{enumerate}[label=\textup{(\roman*)}]
   \item For every \(\lambda\in\Lambda_f\setminus\Gamma\), one has
   \(-\frac{\lambda}{nm}\in\mathcal R_f\).
   \item If \(\Lambda_f\setminus\Gamma\neq\emptyset\) and
   \(\lambda_1=\min(\Lambda_f\setminus\Gamma)\), then
   \(-(\sigma+1)\in\mathcal R_f\) for every
   \(\sigma\in\mathcal S_f\) such that
   \(\sigma<\frac{\lambda_1}{nm}\).
\end{enumerate}
\end{theorem}



 
%



\section{Recovering the first differential value from the Bernstein--Sato roots}

For \(j\in J\), we call \(-s_j\in\mathcal Y_\Gamma\setminus\mathcal T_\Gamma\) the \emph{non-jumped} member of the Cassou--Nogu\`es pair and \(-s_j-1\in\mathcal Q_\Gamma\setminus\mathcal T_\Gamma\) the \emph{jumped} member (see Lemma \ref{lem:combinatorics-root-sets}). The pairs are ordered by the integer \(j\), or equivalently by \(\lambda(j)=j+n+m\).
Analogously,
for \(j\in J\), we call \(-s_j\in \mathcal R_f \cap (\mathcal Y_\Gamma\setminus\mathcal T_\Gamma) \)
a \emph{non-jumped root} of $f$  and \(-s_j-1\in\mathcal R_f \cap (\mathcal Q_\Gamma\setminus\mathcal T_\Gamma)\) a \emph{jumped root}. 
The set 
\begin{equation}
   J(f):=\{j\in J:-s_j\in\mathcal R_f\}.
   \label{eq:J-f-roots}
\end{equation}
indexes injectively the non-jumped Bernstein--Sato roots  of $f$, which are
\begin{equation}
\mathcal R_f \cap (\mathcal Y_\Gamma\setminus\mathcal T_\Gamma) = \{ -s_j: j \in J(f) \}.
   \label{eq:non-jumped-roots}
\end{equation}

Notice that, by Theorem~\ref{thm:cassou-nogues} and Lemma~\ref{lem:combinatorics-root-sets}, in the genus-one case for every \(j\in J\), exactly one of the two rational numbers
\[
   -s_j,
   \qquad -s_j-1
\]
belongs to \(\mathcal R_f\). In particular, the occurrence of the jumped root \(-s_j-1\) excludes the occurrence of the non-jumped root \(-s_j\), and conversely.


\begin{theorem}
\label{MainThm}
\label{thm:main}
Let \(f\) define a plane branch with one Puiseux-pair \(m/n\), \(\gcd(n,m)=1\), \(1<n<m\), and set \(\Gamma=\langle n,m\rangle\). Then \(\mathcal R_f\) determines the first non-topological differential value \(\lambda_1:=\min(\Lambda_f\setminus\Gamma)\), with the convention \(\lambda_1:=\infty\) in the quasihomogeneous case.

More precisely, consider \(J(f)\) defined in \eqref{eq:J-f-roots}.
%
Then
\begin{equation}
   \lambda_1=
   \begin{cases}
      n+m+\min J(f),&\text{if }J(f)\neq\emptyset,\\[2mm]
      \infty,&\text{if }J(f)=\emptyset.
   \end{cases}
   \label{eq:lambda-from-roots}
\end{equation}

Equivalently, the first non-jumped root of \(\mathcal R_f\) in \(\mathcal Y_\Gamma\setminus\mathcal T_\Gamma\), if it exists, is \(-\frac{\lambda_1}{nm}\); if no such root exists, we set \(j_1=\infty\). For every \(\lambda<\lambda_1\) with \(-\frac{\lambda}{nm}\in\mathcal Y_\Gamma\setminus\mathcal T_\Gamma\), the jumped root \(-\frac{\lambda}{nm}-1\) belongs to \(\mathcal R_f\).
\end{theorem}

From this result we obtain the following straightforward consequence:

\begin{corollary}
\label{cor:zariski-from-bs}
The reduced roots \(\mathcal R_f\) of the Bernstein--Sato polynomial determine the Zariski invariant \(j_1(\mathcal R_f)\) of  the genus-one branch $\gamma: f=0$. If 
\(J(f):= \{j\in J:-s_j\in\mathcal R_f\} \neq \emptyset\), then
\[
   j_1(\mathcal R_f) = \min J(f).
\]
Moreover, \(J(f)= \emptyset\) if and only if \( \gamma \) is analytically equivalent to the quasihomogeneous branch.
\end{corollary}

The proof of Theorem~\ref{thm:main} will use the following sharpening of Theorem~\ref{thm:senovilla}. Its point is that the vanishing of all residues before \(j_1\) is not restricted to  \((1,1) \in M\).

\begin{lemma}
\label{lemma-senovilla-improved}
Let \(f\) be as in \eqref{eq:normal-form}. For \(j\in J\), the following are equivalent:
\begin{enumerate}[label=\textup{(\roman*)}]
   \item \(\lambda(j)=j+n+m\) is the first non-topological differential value \(\lambda_1=\min(\Lambda_f\setminus\Gamma)\);
   \item \(j\) is the Zariski invariant \(j_1=\min\{i\in J:z_i\neq0\}\);
   \item for every \((a,b)\in M\), one has
   \[
      \operatorname{Res}_{f}(a,b)(s_\ell)=0
      \qquad\text{for all }\; \ell<j,
      \quad \ell\in J,
   \]
   and
   \[
      \operatorname{Res}_{f}(1,1)(s_j)\neq0.
   \]
\end{enumerate}
\end{lemma}

\begin{proof}
The equivalence between \textup{(i)} and \textup{(ii)} is Theorem~\ref{thm:senovilla}. The implication \textup{(iii)}\(\Rightarrow\)\textup{(ii)} also follows from Theorem~\ref{thm:senovilla}, because \textup{(iii)} contains the corresponding condition for \((a,b)=(1,1)\).

It remains to prove \textup{(ii)}\(\Rightarrow\)\textup{(iii)}. Assume \(j=j_1\). Then \(z_h=0\) for all \(h<j_1\), and \(z_{j_1}\neq0\). Fix \(\ell<j_1\), \(\ell\in J\), and \((a,b)\in M\). In the residue formula \eqref{eq:cassou-residue}, the relevant integer is
\[
   k_\ell(a,b)=\ell+n(1-a)+m(1-b)\leq \ell<j_1,
\]
because \(a,b\geq1\). If \(k_\ell(a,b)<0\), no term appears in the sum. If \(k_\ell(a,b)=0\), then
\[
   \ell=n(a-1)+m(b-1)\in\Gamma,
\]
which is impossible because \(\ell\in J\subset\Gaps(\Gamma)\) by Lemma~\ref{lem:combinatorics-root-sets}. Hence either the sum is empty or \(0<k_\ell(a,b)<j_1\).

In the latter case, every nonzero solution of
\[
   \sum_{h\in J}\delta_hh=k_\ell(a,b)
\]
uses only indices \(h<j_1\). Therefore every monomial
\(\prod_h z_h^{\delta_h}\) appearing in the residue contains at least one factor \(z_h\) with \(h<j_1\), and hence vanishes. Thus
\[
   \operatorname{Res}_{f}(a,b)(s_\ell)=0
\]
for all \(\ell<j_1\) and all \((a,b)\in M\).

Finally, for \((a,b)=(1,1)\) and \(j=j_1\), the integer in \eqref{eq:k-j-ab} is \(k_{j_1}(1,1)=j_1\). In \eqref{eq:cassou-residue}, the only nonzero contribution comes from the solution \(\delta_{j_1}=1\), \(\delta_h=0\) for \(h\neq j_1\). Hence
\[
   \operatorname{Res}_{f}(1,1)(s_{j_1})
   =-\frac{\Gamma(s_{j_1})^{-1}}{nm}
   \Gamma\!\left(\frac{p_{1,j_1}+1}{m}\right)
   \Gamma\!\left(\frac{p_{2,j_1}+1}{n}\right)z_{j_1},
\]
which is nonzero because the Euler gamma function has no zeros and \(z_{j_1}\neq0\).
\end{proof}

\begin{proof}[Proof of Theorem~\ref{thm:main}]
Assume first that \(f\) is not quasihomogeneous, so \(j_1<\infty\). By Lemma~\ref{lemma-senovilla-improved}, for every \(j<j_1\), all residues \(\operatorname{Res}_f(a,b)(s_j)\), \((a,b)\in M\), vanish. Therefore Theorem~\ref{thm:cassou-nogues} gives
\[
   -s_j\notin\mathcal R_f,
   \qquad
   -s_j-1\in\mathcal R_f
   \qquad (j<j_1).
\]
Again by Lemma~\ref{lemma-senovilla-improved},
\[
   \operatorname{Res}_f(1,1)(s_{j_1})\neq0,
\]
and hence Theorem~\ref{thm:cassou-nogues} gives \(-s_{j_1}\in\mathcal R_f\). The same result ensures that these alternatives cannot occur simultaneously within a pair. It follows that
\[
   \min\{j\in J:-s_j\in\mathcal R_f\}=j_1,
\]
and therefore \(\lambda_1=j_1+n+m\) is recovered from \(\mathcal R_f\) by \eqref{eq:lambda-from-roots}.

If \(f\) is quasihomogeneous, then \(z_j=0\) for every \(j\in J\). In the residue formula \eqref{eq:cassou-residue}, every residue \(\operatorname{Res}_f(a,b)(s_j)\) vanishes; hence Theorem~\ref{thm:cassou-nogues} gives \(-s_j-1\in\mathcal R_f\) and \(-s_j\notin\mathcal R_f\) for all \(j\in J\). Thus \(J(f)=\emptyset\), and \eqref{eq:lambda-from-roots} gives \(\lambda_1=\infty\), as required.
\end{proof}




\begin{remark}
Our recovery formula in Theorem~\ref{thm:main} indexes the variable Bernstein--Sato root pairs by
the relevant part of the Yano set
whereas Guimar\~aes--Hernandes use spectral numbers and differential values
\cite[Theorem~4.1, Lemma~4.2 and Proposition~4.4]{guimaraes2025bernstein}.
Their results imply our recovery statement once the topological roots
are identified and the simultaneous occurrence of the jumped and
non-jumped members of each pair is excluded. Thus the recovery
assertion itself is a consequence of their work. Our contribution
here is its explicit formulation through the Yano set and the
alternative residue proof.
\end{remark}

\section{Generic Bernstein--Sato roots at fixed Zariski invariant}

Let \(\mathcal M_\Gamma\) denote the Zariski moduli space of analytic classes of plane
branches with value semigroup \(\Gamma=\langle n,m\rangle\) \cite{zariski1973modules}. For
\(j_1\in J\), let \(\mathcal M_{\Gamma,j_1}\) denote the stratum
consisting of (classes of) branches with Zariski invariant \(j_1\). The resulting
stratification by constant Zariski invariant was studied by Peraire
\cite{peraire1998moduli}.

Theorem~\ref{thm:main} implies that the stratification by constant 
Bernstein--Sato root sets refines the stratification by the Zariski
invariant. Indeed, if \(\mathcal M_{\Gamma,\mathcal R}\) denotes the
locus of (analytic classes of) branches 
whose reduced Bernstein--Sato root set is \(\mathcal R\), then
\[
   \mathcal M_{\Gamma,\mathcal R}
   \subseteq
   \mathcal M_{\Gamma,j_1(\mathcal R)},
\]
where \(j_1(\mathcal R)\) denotes the Zariski invariant recovered from
\(\mathcal R\) by Theorem~\ref{thm:main} (see Corollary~\ref{cor:zariski-from-bs}).

It is therefore natural to ask whether the reduced Bernstein--Sato
root set is constant on a nonempty Zariski-open subset of each stratum
\(\mathcal M_{\Gamma,j_1}\), and whether this generic root set can be
described, in the spirit of Yano's conjecture, in terms of \(\Gamma\)
and \(j_1\). More precisely, for a generic branch in
\(\mathcal M_{\Gamma,j_1}\), does every Cassou--Nogu\`es pair indexed
by \(j\in J\) with \(j\geq j_1\) contribute the non-jumped root
\(-s_j\), rather than the jumped root \(-s_j-1\)? The following
result answers this question.

\begin{theorem}
\label{thm:generic-roots}
Let \(j_1\in J\) and put \(\lambda_1=j_1+n+m\). For a generic branch in \(\mathcal M_{\Gamma,j_1}\), the reduced Bernstein--Sato root set is the disjoint union
\begin{equation}
   \mathcal R(\Gamma,j_1)
   =\mathcal R^{\mathrm{top}}_\Gamma
    \sqcup\mathcal R^{<\lambda_1}_\Gamma
    \sqcup\mathcal R^{\geq\lambda_1}_\Gamma,
   \label{eq:generic-root-union}
\end{equation}
where
\begin{align*}
   \mathcal R^{\mathrm{top}}_\Gamma
      &:=\mathcal T_\Gamma,\\
   \mathcal R^{<\lambda_1}_\Gamma
      &:=\left\{
          -\frac{\lambda(j)}{nm}-1:
          j\in J,
          j<j_1
        \right\}
        =\left\{
          -\frac{\lambda}{nm}-1:
          \lambda<\lambda_1,
          -\frac{\lambda}{nm}\in\mathcal Y_\Gamma\setminus\mathcal T_\Gamma
        \right\},\\
   \mathcal R^{\geq\lambda_1}_\Gamma
      &:=\left\{
          -\frac{\lambda(j)}{nm}:
          j\in J,
          j\geq j_1
        \right\}
        =\left\{
          -\frac{\lambda}{nm}:
          \lambda\geq\lambda_1,
          -\frac{\lambda}{nm}\in\mathcal Y_\Gamma\setminus\mathcal T_\Gamma
        \right\}.
\end{align*}
Equivalently, in a generic \(\mu\)-constant and \(j_1\)-constant deformation, the pairs indexed by \(j<j_1\) jump, and all pairs indexed by \(j\geq j_1\) do not jump.
\end{theorem}

\begin{proof}
Let
\[
   A_{\Gamma,j_1}:=
   \left\{\mathbf z=(z_j)_{j\in J}\in\C^J:
      z_h=0\ \text{for }h<j_1,\quad z_{j_1}\neq0
   \right\}.
\]
This is the irreducible coefficient space
\[
   A_{\Gamma,j_1}\simeq
   \C^*\times\C^{\#\{j\in J:j>j_1\}}
\]
of normal-form equations
\begin{equation}
   f_{\mathbf z}
   =x^m+y^n+
     \sum_{\substack{j\in J\\j\geq j_1}}
       z_jx^{p_{1,j}}y^{p_{2,j}}.
   \label{eq:j1-stratum-normal-form}
\end{equation}
Let
\[
   \pi:A_{\Gamma,j_1}\longrightarrow\mathcal M_{\Gamma,j_1}
\]
be the quotient map from the normal-form construction of the moduli
stratum \cite[\S4]{peraire1998moduli}; it sends an equation to its
analytic equivalence class. 
For every branch in this stratum, Theorem~\ref{thm:main} implies that \(-s_j-1\in\mathcal R_f\) for all \(j<j_1\), while \(-s_{j_1}\in\mathcal R_f\).

It remains to construct a dense open subset of the moduli stratum on
which \(-s_{j_2}\) occurs for every \(j_2>j_1\). Fix such an index and
consider the residue
\[
   R_{j_2}(\mathbf z):=\operatorname{Res}_{f_{\mathbf z}}(1,1)(s_{j_2}).
\]
By \eqref{eq:cassou-residue}, this is a polynomial in the parameters \(z_j\). The monomial \(z_{j_2}\) appears from the solution
\[
   \delta_{j_2}=1,
   \qquad \delta_h=0\quad(h\neq j_2)
\]
of \(\sum_h\delta_hh=j_2\). In fact, this is the unique summand
involving \(z_{j_2}\). Indeed, every \(h\in J\) is positive, so any
solution with \(\delta_{j_2}>0\) satisfies
\[
   j_2=\sum_{h\in J}\delta_hh
      \geq \delta_{j_2}j_2.
\]
It follows that \(\delta_{j_2}=1\), and equality then forces
\(\delta_h=0\) for every \(h\neq j_2\). Thus no other summand can
contribute to the coefficient of \(z_{j_2}\). Consequently,
\[
   R_{j_2}(\mathbf z)
      =C_{j_2}z_{j_2}
       +P_{j_2}\bigl((z_h)_{j_1\leq h<j_2}\bigr),
\]
where \(P_{j_2}\) is independent of \(z_{j_2}\) and
\[
   C_{j_2}=-\frac{\Gamma(s_{j_2})^{-1}}{nm}
   \Gamma\!\left(\frac{p_{1,j_2}+1}{m}\right)
   \Gamma\!\left(\frac{p_{2,j_2}+1}{n}\right).
\]
Therefore \(R_{j_2}(\mathbf z)\) is not the zero
polynomial. Since \(A_{\Gamma,j_1}\) is irreducible, the principal open
\[
   D(R_{j_2}):=
   \left\{\mathbf z\in A_{\Gamma,j_1}:
      R_{j_2}(\mathbf z)\neq0
   \right\}
\]
is nonempty and dense in \(A_{\Gamma,j_1}\).

To descend this statement to the quotient, consider the 
root-occurrence locus
\[
   G_{j_2}:=
   \left\{\mathbf z\in A_{\Gamma,j_1}:
      -s_{j_2}\in\mathcal R_{f_{\mathbf z}}
   \right\}.
\]
By Theorem~\ref{thm:cassou-nogues},
\[
   G_{j_2}
   =\bigcup_{(a,b)\in M}
      \left\{\mathbf z\in A_{\Gamma,j_1}:
         \operatorname{Res}_{f_{\mathbf z}}(a,b)(s_{j_2})\neq0
      \right\}.
\]
Thus \(G_{j_2}\) is Zariski open, and it is dense because it contains
\(D(R_{j_2})\). Moreover, \(G_{j_2}\) is saturated with respect to
\(\pi\) (equals the inverse image of its own image by \(\pi\)), since the Bernstein--Sato polynomial is invariant under
analytic equivalence. Hence
\[
   G:=\bigcap_{\substack{j_2\in J\\j_2>j_1}}G_{j_2}
\]
is a dense Zariski-open saturated subset of \(A_{\Gamma,j_1}\).

Set \(U:=\pi(G)\). Since \(G\) is saturated,
\(\pi^{-1}(U)=G\); consequently, \(U\) is open in the quotient topology
of \(\mathcal M_{\Gamma,j_1}\). It is also dense. Indeed, if
\(V\subseteq\mathcal M_{\Gamma,j_1}\) is nonempty and open, then
\(\pi^{-1}(V)\) is a nonempty open subset of \(A_{\Gamma,j_1}\), so it
meets \(G\); hence \(V\cap U\neq\emptyset\).

On \(U\), all roots indexed by \(j\geq j_1\) are non-jumped, while all
roots indexed by \(j<j_1\) are jumped.
Finally, the Cassou--Nogu\`es decomposition
\eqref{eq:cassou-root-decomposition} and Lemma ~\ref{lem:combinatorics-root-sets}
 show that the displayed union in \eqref{eq:generic-root-union} accounts for all
reduced roots.
\end{proof}

\begin{remark}
The quasihomogeneous case corresponds to  \(j_1=\lambda_1=\infty\). Then every pair is jumped, and Theorem~\ref{thm:cassou-nogues} recovers
\[
   \mathcal R_f=\mathcal T_\Gamma
   \sqcup\left\{-\frac{\lambda(j)}{nm}-1:j\in J\right\}
   =\mathcal Q_\Gamma.
\]
At the opposite extreme, the generic stratum of the whole topological class is obtained by taking the smallest possible Zariski invariant; in that case the generic Bernstein--Sato roots are \(\mathcal Y_\Gamma\).
\end{remark}

\begin{remark}
The open set \(U\) constructed in the proof of
Theorem~\ref{thm:generic-roots} is exactly
\(\mathcal M_{\Gamma,\mathcal R(\Gamma,j_1)}\). Hence this locus is
open and dense in \(\mathcal M_{\Gamma,j_1}\), and it is locally
closed in \(\mathcal M_\Gamma\) by
\cite[Proposition~4.5]{peraire1998moduli}. With the constructible
complex structures fixed above, it is an open constructible subspace
of \(\mathcal M_{\Gamma,j_1}\). In the coefficient space,
\[
   \operatorname{codim}_{\C^J}\overline{A_{\Gamma,j_1}}
   =\#\{j\in J:j<j_1\}
   =\#\mathcal R^{<\lambda_1}_\Gamma.
\]
Thus the number of generically jumped Bernstein--Sato roots equals this coefficient-space
codimension.
\end{remark}

\begin{remark}\label{rem:saito-specialization}
Theorem~\ref{thm:generic-roots} gives an explicit description of the generic
roots in every stratum \(\mathcal M_{\Gamma,j_1}\).
Its coefficient-space statement can also be obtained from Saito's
deformation formulas \cite{saito2022semiweighted}. Our residue argument makes the
description in terms of \((\Gamma,j_1)\) explicit and establishes the
corresponding dense open locus in the moduli quotient.
\end{remark}

\section{Spectral--exponent discrepancy and the dimension of generic-root loci}
\label{sec:spectral-exponent-discrepancy}

In this section, we study the spectral--exponent discrepancy, an analytical invariant that counts the number of \(b\)-exponents whose opposites are spectral numbers of an isolated singularity $f$. We introduce the spectral--exponent discrepancy in the setting where the set of \(b\)-exponents coincides with the set of reduced Bernstein--Sato roots, \(\mathcal B_f=\mathcal R_f\), as is the case for genus-one plane branches.
In this case, we first express the spectral--exponent discrepancy as a count of non-jumped roots, and then use Theorem~\ref{thm:generic-roots} to compute its generic value. 

Let \(\mathcal S_f\) be the spectrum, normalized to lie in \((-1,1)\), of
the branch defined by $f$  and let
\[
   \mathcal E_f:=\{-(\rho+1):\rho\in\mathcal B_f\} = -(\mathcal B_f+1) = -(\mathcal R_f+1) 
\]
be the set of \textit{exponents} associated with the \(b\)-exponents and hence with the reduced Bernstein--Sato roots of $f$ in our setting. 
The spectrum of a plane branch is determined by its Puiseux
pairs (see Saito's formula \cite[Theorem~1.5]{saito2000exponents}) and is hence a topological invariant.
For any branch with semigroup $\Gamma=\langle n,m\rangle$,  Saito's formula for one-Puiseux-pair,
rewritten in terms of semigroup gaps (see also \cite[Lemma 4.2]{guimaraes2025bernstein}), gives
\[
\mathcal S_f= \mathcal S_\Gamma =
\left\{\pm\frac{\delta}{nm}:
0 <\delta\in \Gaps(\Gamma)
\right\}.
\]

Following Saito and Hertling--Stahlke
\cite{saito1989brieskorn,hertling1999bernstein}, consider the
spectral--exponent discrepancy (originally defined as the dimension of the
quotient of the saturated Brieskorn lattice by the Brieskorn lattice):
\[
   D_f:=
   \sum_{\sigma\in\mathcal S_f}\sigma
   -\sum_{\epsilon\in\mathcal E_f}\epsilon.
\]
It satisfies \cite[Proposition~3.5]{hertling1999bernstein}
\(
\mu_f-\tau_f=\#(\Lambda_f\setminus\Gamma_f)\leq D_f.
\)
To interpret this difference of sums in the genus-one case, we write each
sum using the sets introduced earlier.

For a quasihomogeneous branch defined by $f_0=x^m+y^n$, it is well knwown that the roots of the reduced
Bernstein--Sato polynomial satisfy
$\mathcal R_{f_0}=-(\mathcal S_{f_0}+1)$; see
\cite[Corollary to Theorem 14.14
]{arnold1988singularitiesII},
applying the corresponding shift in the spectrum convention. Hence the spectral--exponent discrepancy is zero for any quasihomogeneous branch, since its exponents coincide with the spectrum, 
and we obtain
\[ \mathcal S_{\Gamma}=-(\mathcal R_{f_0}+1)=-(\mathcal Q_{\Gamma}+1).\]
Since \(\mathcal Q_\Gamma\) is
the disjoint union of the topological roots $\mathcal T_\Gamma$ and the roots
\(-s_j-1\), \(j\in J\) (see Lemma~\ref{lem:combinatorics-root-sets}), we obtain
\[
   \sum_{\sigma\in\mathcal S_f}\sigma
   =-\sum_{\rho\in\mathcal T_\Gamma}(\rho+1)
     +\sum_{j\in J}s_j.
\]

For the exponents sum, recall from \eqref{eq:J-f-roots} and  \eqref{eq:non-jumped-roots} that
\(J(f)=\{j\in J:-s_j\in\mathcal R_f\}\) indexes injectively the non-jumped Bernstein--Sato roots 
\(
\mathcal R_f \cap (\mathcal Y_\Gamma\setminus\mathcal T_\Gamma) = \{ -s_j: j \in J(f) \}.
\)
Thus
\[
   \sum_{\epsilon\in\mathcal E_f}\epsilon
   =-\sum_{\rho\in\mathcal T_\Gamma}(\rho+1)
     +\sum_{j\in J(f)}(s_j-1)
     +\sum_{j\in J\setminus J(f)}s_j.
\]
Subtracting, the contributions from \(\mathcal T_\Gamma\) cancel,
as do those indexed by \(J\setminus J(f)\). Each index in \(J(f)\)
contributes exactly \(s_j-(s_j-1)=1\), so
\begin{equation}
   D_f=\#J(f)
   =\#\bigl(\mathcal R_f\cap
      (\mathcal Y_\Gamma\setminus\mathcal T_\Gamma)\bigr).
   \label{eq:discrepancy-counts-nonjumped}
\end{equation}
Thus \(D_f\) counts the non-topological, non-jumped Bernstein--Sato
roots of $f$: each such root contributes one to the
discrepancy, while topological roots and jumped roots contribute zero.

For \(j_1\in J\), put \(\lambda_1=j_1+n+m\) and define
\begin{equation}
   D_{\mathrm{gen}}(\Gamma,j_1)
   :=
   \#\bigl(
      \mathcal R(\Gamma,j_1)
      \cap(\mathcal Y_\Gamma\setminus\mathcal T_\Gamma)
   \bigr).
   \label{eq:D-generic-definition}
\end{equation}
By Theorem~\ref{thm:generic-roots}, the non-jumped roots for a
generic branch in \(\mathcal M_{\Gamma,j_1}\) are precisely those
indexed by \(j\geq j_1\). Consequently,
\begin{equation}
\begin{aligned}
   D_{\mathrm{gen}}(\Gamma,j_1)
   &=\#\mathcal R^{\geq\lambda_1}_\Gamma
     =\#\{j\in J:j\geq j_1\}\\
   &=\#\{\ell\in\Gaps(\Gamma):\ell\geq\lambda_1\}.
\end{aligned}
\label{eq:D-generic-count}
\end{equation}
The last equality uses the bijection \eqref{eq:J-large-gaps} of
Lemma~\ref{lem:combinatorics-root-sets}: the inequality \(j\geq j_1\)
is equivalent to \(\lambda(j)=j+n+m\geq\lambda_1\).
By \eqref{eq:discrepancy-counts-nonjumped},
\(D_{\mathrm{gen}}(\Gamma,j_1)\) is therefore the value of \(D_f\)
on the generic root locus in \(\mathcal M_{\Gamma,j_1}\).

\begin{theorem}
\label{Th_Spectral-BS-Discrepancy}
\label{thm:discrepancy-bound}
Let \(j_1\in J\), let \(\lambda_1=j_1+n+m\), and let \(f\) define 
a branch in \(\mathcal M_{\Gamma,j_1}\). Then
\begin{equation}
   D_f
   \leq
   D_{\mathrm{gen}}(\Gamma,j_1).
   \label{eq:D-bound-gaps}
\end{equation}
Moreover,
\begin{equation}
   D_{\mathrm{gen}}(\Gamma,j_1)
   =
   \sum_{a=0}^{\left\lfloor\lambda_1/n\right\rfloor}
      \left(
         \left\lfloor\frac{\lambda_1-an}{m}\right\rfloor+1
      \right)
   -\lambda_1
   +\frac{(n-1)(m-1)}{2}.
   \label{eq:D-bound-formula}
\end{equation}
The inequality \eqref{eq:D-bound-gaps} is sharp: equality holds for a generic branch in \(\mathcal M_{\Gamma,j_1}\).
\end{theorem}

\begin{proof}
By \eqref{eq:discrepancy-counts-nonjumped}, \(D_f\) is the number of roots of \(\mathcal R_f\) lying in \(\mathcal Y_\Gamma\setminus\mathcal T_\Gamma\). By Theorem~\ref{thm:main}, no element \(-\frac{\lambda(j)}{nm}\) with \(j<j_1\) can be such a root for a branch in \(\mathcal M_{\Gamma,j_1}\). Hence, by \eqref{eq:D-generic-count},
\[
   D_f
   \leq
   \#\left\{j\in J:j\geq j_1\right\}
   =D_{\mathrm{gen}}(\Gamma,j_1).
\]
This proves \eqref{eq:D-bound-gaps}. Theorem~\ref{thm:generic-roots}
gives equality on the generic open subset of
\(\mathcal M_{\Gamma,j_1}\), because there all pairs with
\(j\geq j_1\) are non-jumped.

It remains to compute the cardinality of the set of gaps at least \(\lambda_1\). Since \(\#\Gaps(\Gamma)=\frac{(n-1)(m-1)}{2}\), we have
\begin{equation}
\begin{aligned}
   \#\{\ell\in\Gaps(\Gamma):\ell\geq\lambda_1\}
   &=\frac{(n-1)(m-1)}{2}
     -\#\{\ell\in\Gaps(\Gamma):\ell<\lambda_1\}\\
   &=\frac{(n-1)(m-1)}{2}
     -\lambda_1
     +\#\{\ell\in\Gamma:\ell<\lambda_1\}.
\end{aligned}
\label{eq:gaps-count-step}
\end{equation}
Because \(\lambda_1\) is a gap of \(\Gamma\), the conditions \(\ell<\lambda_1\) and \(\ell\leq\lambda_1\) are equivalent for \(\ell\in\Gamma\). Thus \(\#\{\ell\in\Gamma:\ell<\lambda_1\}\) is the number of pairs \((a,b)\in\Z_{\geq0}^2\) such that
\[
   an+bm\leq\lambda_1.
\]
There is no overcounting in this range: since \(\lambda_1\leq nm-n-m<nm\), two representations \(an+bm=a'n+b'm< nm\) would force \(m\mid a-a'\), and our range has \(0\leq a,a'<m\). Hence \(a=a'\) and then \(b=b'\).
For each fixed \(a\) with \(0\leq a\leq\lfloor\lambda_1/n\rfloor\), the number of possible \(b\)'s is
\[
   \left\lfloor\frac{\lambda_1-an}{m}\right\rfloor+1.
\]
Substituting this count into \eqref{eq:gaps-count-step}, and using
\eqref{eq:D-generic-count}, gives \eqref{eq:D-bound-formula}.
\end{proof}

Together with the Hertling--Stahlke inequality
\cite[Proposition~3.5]{hertling1999bernstein}, Theorem~\ref{thm:discrepancy-bound}
gives
\begin{equation} 
\label{ineqDiscrepancy}
   \mu-\tau_f
   =\#(\Lambda_f\setminus\Gamma)
   \leq D_f
   \leq D_{\mathrm{gen}}(\Gamma,j_1).
\end{equation}
The second inequality is an equality on the generic root locus, and the
last term is given explicitly by \eqref{eq:D-bound-formula}.


\begin{remark}
The numerical upper bound for the discrepancy is already implicit in
\cite[(31)]{guimaraes2025bernstein}. We identify its generic attainment
and express it as a semigroup-gap count and a closed floor-sum (Theorem \ref{thm:discrepancy-bound}).
\end{remark}

The following result shows how \(D_{\mathrm{gen}}(\Gamma,j_1)\) determines
the dimensions of both the stratum \(\mathcal M_{\Gamma,j_1}\) and its
Bernstein--Sato generic-root locus \(\mathcal M_{\Gamma,\mathcal R(\Gamma,j_1)}\).

\begin{corollary}
\label{cor:generic-root-locus-dimension}
Let \(j_1\in J\), and set
\[
   \tau_{\min}(\Gamma,j_1)
   :=
   \min\{\tau_\gamma:[\gamma]\in\mathcal M_{\Gamma,j_1}\}.
\]
Then
\begin{equation}
   \dim\mathcal M_{\Gamma,\mathcal R(\Gamma,j_1)}
   =\dim\mathcal M_{\Gamma,j_1}
   = D_{\mathrm{gen}}(\Gamma,j_1)
   -\bigl(\mu-\tau_{\min}(\Gamma,j_1)\bigr).
   \label{eq:generic-root-locus-dimension}
\end{equation}
\end{corollary}

\begin{proof}
Put
\[
   X:=\mathcal M_{\Gamma,j_1},\qquad
   B:=\mathcal M_{\Gamma,\mathcal R(\Gamma,j_1)},\qquad
   T:=\{[\gamma]\in X:\tau_\gamma=\tau_{\min}(\Gamma,j_1)\}.
\]
By Theorem~\ref{thm:generic-roots}, \(B\) is open and dense in
\(X\). By \cite[Proposition~4.7]{peraire1998moduli}, \(T\) is also
open and dense in \(X\). These are constructible subspaces of the
complex structure fixed above. Since the dimension of a
finite-dimensional complex constructible space is unchanged upon
passage to a dense open subspace,
\[
   \dim B=\dim X=\dim T.
\]
Peraire's invariant \(s\) is \(j_1\) in the present notation, and
her integer \(N_s\) is \(\#\{j\in J:j>j_1\}\). Her dimension
theorem therefore gives
\[
   \dim T=N_s+1-\mu+\tau_{\min}(\Gamma,j_1)
\]
\cite[Theorem~4.12]{peraire1998moduli}. Finally,
\[
   N_s+1
   =\#\{j\in J:j\geq j_1\}
   =D_{\mathrm{gen}}(\Gamma,j_1)
\]
by \eqref{eq:D-generic-count}, which proves
\eqref{eq:generic-root-locus-dimension}.
\end{proof}

Peraire
\cite[\S~2 and \S~3]{peraire1998moduli} gives an algorithmic construction of a semimodule of differential values \( \Lambda_{\mathrm{gen}}(\Gamma , j_1)\) attained by generic branches in \(\mathcal M_{\Gamma,j_1}\). Then she establishes
\[ \mu - \tau_{\min}(\Gamma,j_1) =\#\bigl(\Lambda_{\mathrm{gen}}(\Gamma,j_1)\setminus\Gamma\bigr).\]

\begin{corollary}
With \(\lambda_1=j_1+n+m\),
\[
   \dim\mathcal M_{\Gamma,\mathcal R(\Gamma,j_1)}
   =\dim\mathcal M_{\Gamma,j_1}
   =\#\{\ell\in\Gaps(\Gamma):
          \ell\geq\lambda_1,\ \ell\notin\Lambda_{\mathrm{gen}}(\Gamma,j_1)\}.
\]
\end{corollary}

\begin{proof}
By Corollary~\ref{cor:generic-root-locus-dimension} and
\eqref{eq:D-generic-count},
\[
\begin{aligned}
   \dim\mathcal M_{\Gamma,\mathcal R(\Gamma,j_1)}
   &=D_{\mathrm{gen}}(\Gamma,j_1)
     -\bigl(\mu-\tau_{\min}(\Gamma,j_1)\bigr)\\
   &=\#\{\ell\in\Gaps(\Gamma):\ell\geq\lambda_1\}
     -\#\bigl(\Lambda_{\mathrm{gen}}(\Gamma,j_1)\setminus\Gamma\bigr).
\end{aligned}
\]
Since \(\Lambda_{\mathrm{gen}}(\Gamma,j_1)\setminus\Gamma\)
consists of gaps at least \(\lambda_1\), the difference of these
cardinalities is the stated count.
\end{proof}

\section{Optimality and examples}

Theorem~\ref{thm:main} is optimal is optimal in the sense that higher non-topological differential values cannot be recovered from the Bernstein-Sato roots in general, as next example shows.


\begin{example}
\label{ex:optimality-5-12}
Consider the topological class \(\Gamma=\langle5,12\rangle\). Let
\[
   f_1=y^5-x^{12}+x^5y^3,
   \qquad
   f_2=y^5-x^{12}+x^5y^3+\frac{3}{10}x^{10}y.
\]
Both branches have the same Zariski invariant and the same first non-topological differential value \(\lambda_1=18\). A direct computation using the Cassou--Nogu\`es criterion gives
\[
   \mathcal R_{f_1}=\mathcal R_{f_2}=\mathcal Y_{\langle5,12\rangle}.
\]
However, their semimodules of differential values are different:
\[
   \Lambda_{f_1}=\langle5,12,18,31\rangle,
   \qquad
   \Lambda_{f_2}=\langle5,12,18\rangle.
\]
Thus the two branches have the same reduced Bernstein--Sato root set, while the higher non-topological differential value \(31\) occurs for \(f_1\) and not for \(f_2\). Consequently, no statement stronger than Theorem~\ref{thm:main} can hold in general.
\end{example}

\begin{remark}
\label{rem:tjurina-root-strata-incomparable}
Neither the Tjurina stratification nor the reduced Bernstein--Sato
stratification refines the other. Indeed,
\[
   \#(\Lambda_{f_1}\setminus\Gamma)=7,
   \qquad
   \#(\Lambda_{f_2}\setminus\Gamma)=6.
\]
Since \(\mu=44\) and
\(\mu-\tau_f=\#(\Lambda_f\setminus\Gamma)\), the two branches in
Example~\ref{ex:optimality-5-12} have the same reduced root set but
Tjurina numbers \(37\) and \(38\), respectively. Conversely,
Peraire exhibits two branches with semigroup \(\langle5,11\rangle\)
and Tjurina number \(34\), but with Zariski invariants \(1\) and
\(2\) \cite[Remark~4.13]{peraire1998moduli}. Their reduced root sets
are different by Corollary~\ref{cor:zariski-from-bs}.
\end{remark}

For \(\Gamma=\langle5,12\rangle\), the set \(J\) is
\[
   J=\{1,2,4,6,9,11,14,16,21,26\},
\]
and the corresponding values \(\lambda(j)=j+17\) are
\[
   18,19,21,23,26,28,31,33,38,43.
\]
These are precisely the gaps of \(\langle5,12\rangle\) that are at least \(17=n+m\). Hence, if a branch in this class has Zariski invariant \(j_1\), Theorem~\ref{thm:generic-roots} says that the generic Bernstein--Sato root set in \(\mathcal M_{\Gamma,j_1}\) is obtained from the Yano set by shifting exactly the roots indexed by the elements of \(J\) smaller than \(j_1\).

\begin{example}[Inclusion between  generic-root and generic-differential-values strata]
\label{ex:fixed-j1-14-boundary}
Still in the topological class \(\Gamma=\langle5,12\rangle\), consider
\[
\begin{aligned}
   g_+&=x^{12}+y^5+x^{10}y^2+x^8y^3,\\
   g_1&=x^{12}+y^5+x^{10}y^2+x^9y^3,\\
   g_0&=x^{12}+y^5+x^{10}y^2.
\end{aligned}
\]
For all three branches the first nonzero standing coefficient is the one
indexed by \(j=14\). Thus
\[
   j_1=14,
   \qquad
   \lambda_1=31.
\]
The corresponding Tjurina numbers of \(g_+, g_1, g_0 \) are \(41,42,42\).
Since \(\mu=44\), it
follows that
\[
   \mu-\tau_{g_+}=3,
   \qquad
   \mu-\tau_{g_1}=\mu-\tau_{g_0}=2.
\]
%
The three branches lie in the only two differential-values strata of
\(\mathcal M_{\Gamma,14}\): 
\[   \Lambda_{g_+}=\langle5,12,31,38\rangle,
   \qquad
   \Lambda_{g_1}=\Lambda_{g_0}=\langle5,12,31\rangle.
\]
%
In the normal form
\eqref{eq:j1-stratum-normal-form}, the relevant coefficients for
\(g_+,g_1,g_0\), respectively, are
\[
   (z_{14},z_{16},z_{21})=(1,1,0),
   \qquad (1,0,1),
   \qquad (1,0,0),
\]
with all other \(z_j\) equal to zero.
From Cassou-Noguès criterion we infer that the pair indexed by \(j=16\) contributes the non-jumped Bernstein--Sato root
\(-s_{16}=-11/20\) when \(z_{16}\neq0\), whereas the jumped root
\(-s_{16}-1=-31/20\) occurs when \(z_{16}=0\); the pair indexed by \(j=21\) contributes \(-19/30\) when either
\(z_{16}\neq0\) or \(z_{21}\neq0\), and contributes its translate
\(-49/30\) when both coefficients vanish; the pairs indexed by \(j=14\) and \(j=26\) do not jump. 
All pairs with \(j<14\) are jumped by Theorem~\ref{thm:main}.
Hence,
the only Bernstein--Sato roots that differ among these branches are summarized by
\[
\begin{array}{c|c|c|c|c|c}
 & \Lambda & \text{root from }j=16 & \text{root from }j=21
 & \mu-\tau & D_f\\
\hline
g_+ & \langle5,12,31,38\rangle & -\frac{11}{20} & -\frac{19}{30}
    & 3 & 4\\
g_1 & \langle5,12,31\rangle & -\frac{31}{20} & -\frac{19}{30}
    & 2 & 3\\
g_0 & \langle5,12,31\rangle & -\frac{31}{20} & -\frac{49}{30}
    & 2 & 2
\end{array}
\]
In this \(\mathcal M_{\Gamma,14}\) stratum \(D_{\mathrm{gen}}(\Gamma,14)=4\). Hence the three
rows exhibit, respectively, the sharp upper bound \(D_f=D_{\mathrm{gen}}\),
a strict instance of both inequalities \eqref{ineqDiscrepancy}, and equality in the lower bound
\(\mu-\tau_f\leq D_f\).

This example illustrates the new
generic-root and discrepancy results of this paper, and it complements the previous results
of Guimar\~aes and Hernandes, which already showed that
the differential-values and Bernstein--Sato-root stratifications were not comparable
\cite[Theorem~4.1 and Example~5.1]{guimaraes2025bernstein}.
Throughout \(\mathcal M_{\Gamma,14}\), the \(j=16\) root pair
distinguishes the two differential values substrata:
\(-11/20\in\mathcal R_f\) if and only if \(38\in\Lambda_f\).
Consequently, \(\mathcal R_f\) determines \(\Lambda_f\) throughout
this stratum, even though the generic inequality is strict:
\[
D_{\mathrm{gen}}(\Gamma,14)=4
>3=\mu-\tau_{\mathrm{min}}(\Gamma,14).
\]
To our knowledge, this coexistence of complete recovery of the
differential semimodule with a strict generic inequality has not
been explicitly illustrated in previous works.
Conversely, the comparison of \(g_1\) and \(g_0\) shows that
\(\Lambda_f\) does not determine \(\mathcal R_f\): the
Bernstein--Sato root set can vary within the differential stratum
\(\Lambda_f=\langle5,12,31\rangle\).

In this example the generic-root substratum of  \(\mathcal M_{\Gamma,14}\), i.e. the open subset of branches with Bernstein--Sato roots equal to \(\mathcal R(\Gamma,14)\)), is included in the generic differential-values substratum, i.e. the open subset of branches with differential values equal to \(\Lambda _{\text{gen}}({\Gamma, j_1})\).
\end{example}

\end{document}